\documentclass[a4paper, oneside, 11pt]{article}
\usepackage[english]{babel}
\usepackage{amsfonts}
\usepackage{amssymb}
\usepackage{graphics}
\usepackage{amsrefs}
\usepackage{latexsym}
\begin{document}

\title{{\bf{\Large{A note on a Bayesian nonparametric  estimator of the discovery probability}}}\footnote{{\it AMS (2000) subject classification}. Primary: 60G58. Secondary: 60G09.}}
\author{\textsc {Annalisa Cerquetti}\footnote{Corresponding author, SAPIENZA University of Rome, Via del Castro Laurenziano, 9, 00161 Rome, Italy. E-mail: {\tt annalisa.cerquetti@gmail.com}}\\
\it{\small Department of Methods and Models for Economics, Territory and Finance}\\
  \it{\small Sapienza University of Rome, Italy }}
\newtheorem{teo}{Theorem}
\date{\today}
\maketitle{}

\begin{abstract}
Favaro, Lijoi, and Pr\"unster (2012, {\it Biometrics}, {\bf 68}, 1188--1196) derive a novel Bayesian nonparametric estimator of the probability of detecting at the $(n+m+1)$th observation a species already observed with any given frequency in an enlarged sample of size $n+m$, conditionally on a basic sample of size $n$. Unfortunately the general result under Gibbs priors (Theorem 2), and consequently the explicit result under $(\alpha, \theta)$ Poisson-Dirichlet priors (Proposition 3), appear to be wrong. Here we provide the correct formulas for both the results, obtained by means of a new technique devised in Cerquetti (2013). We verify the correctness of our derivation by an explicit counterproof for the two-parameter Poisson-Dirichlet case. 
\end{abstract}


\section{Introduction}
\label{s:intro}

In species sampling problems, particularly in ecological or biological studies, given a vector $(n_1, \dots, n_j)$,  $\sum_{i=1}^{j} n_i=n$, reporting the multiplicities of $j$ different species observed in order of appearance among the first $n$ observations, interest may lie in estimating the {\it probability} to detect, at step $n+m+1$, a species already observed with any given frequency  both among the species belonging to the basic $n-$sample, or among the {\it new} species eventually arising in the additional $m-$sample which is still to be observed. In  Favaro, Lijoi, and Pr\"unster (2012) a Bayesian nonparametric solution is provided as a generalization to the problem of estimating the {\it discovery probability}, i.e. the probability to discover a {\it new} species not represented in the previous $n+m$ observations, already solved in Lijoi, Mena, and Pr\"unster (2007). See the original article for a discussion and a review of previous proposals in both frequentist and Bayesian approaches. 
As a common {\it prior} assumption the authors assume that the theoretically infinite sequence of unknown relative abundances of the species in the population, $(P_i)_{i\geq 1}$,   follows a random discrete distribution belonging to the {\it Gibbs family}, (Gnedin and Pitman, 2006) i.e. such that, by Kingman's correspondence (Kingman, 1978), the probability to observe a {\it specific} partition of the first $n$ observations among the first $j$ species with multiplicities $(n_1, \dots, n_j)$ (called the exchangeable partition probability function, EPPF) is given by
\begin{equation}
\label{kingm}
p(n_1, \dots, n_j)=\sum_{(l_1, \dots, l_j)} \mathbb{E} \left[  \prod_{i=1}^j P_{l_i}^{n_i} \right]= V_{n,j} \prod_{i=1}^j (1-\alpha)_{n_i-1},
\end{equation}
where $(l_1, \dots, l_j)$ ranges over all ordered $j$-tuples of distinct positive integers, and $(x)_y=x(x+1)\cdots(x+y-1)$ stands for rising factorials. Here the $V=(V_{n,j})$ represent specific Gibbs coefficients identifying a specific Gibbs model belonging to one of the three classes devised by Gnedin and Pitman (2006), and arise as mixtures of three different extreme partitions coefficients for $\alpha \in (-\infty, 0)$, $\alpha=0$ and $\alpha \in (0,1)$. This is equivalent to assume that the theoretically infinite sequence of species {\it labels} $(X_i)_{i \geq 1}$ is exchangeable with almost surely discrete {\it de Finetti measure} representable as $P(\cdot)=\sum_{i=1}^{\infty} P_i \delta_{Y_i}(\cdot)$, for $(P_i)$ any rearrangement of the ranked frequencies $(P_i^{\downarrow})$ satisfying (\ref{kingm}), independent of $(Y_i) \sim$ IID $H(\cdot)$, for $H$ some non atomic probability distribution.

While testing, through an example, the computational advantage provided by the technique based on marginals of multivariate Gibbs distributions devised in Cerquetti (2013), in deriving Bayesian nonparametric estimators in species sampling problems under Gibbs priors,  we find the result in Theorem 2. under general {\it Gibbs priors}, and consequently the particular result in Proposition 3. under {\it two-parameter Poisson-Dirichlet priors}, (Pitman and Yor, 1997) in Favaro et al. (2012) to be affected by a mistake. Despite it is just a subtle problem in a summation limit, it actually produces concrete consequences on calculation of the estimator, as from the $(\alpha, \theta)$ particular case. The aim of this note is to provide potential readers of the Biometrics paper with the correct formulas, in view of possible applications of the Favaro et al. (2012) results under different classes of Gibbs priors and different datasets, or for further derivations of related theoretical results.\footnote{The problem has been highlighted to the authors, and on January 2013 this note has been submitted as a {\it Reader Reaction} to {\it Biometrics}. It has been rejected with the motivation that the numerical results contains in Section 3 of the original paper - an application to genomic data under two-parameter Poisson Dirichlet priors - are not significantly affected by the mistaken formula of the corresponding estimator.}  

We start by noticing that the complexity of the proof of Theorem 1. obtained in Favaro et al. (2012, cfr. Web Appendix pag. 2), and concerning the Bayesian nonparametric estimator of the {\it $[0:k]$-discovery}, namely the probability to observe at step $n+1$ a species already represented $k$ times among the first $n$ observations, may be avoided by encoding the EPPF of a general Gibbs partition, $p(n_1, \dots, n_j)=V_{n,j}\prod_{i=1}^j (1-\alpha)_{n_i-1}$, in terms of the counting vector of blocks/species of different sizes, namely
$$
q(c_1, \dots, c_{n})= V_{n,j} \prod_{k=1}^{n} [(1-\alpha)_{k-1}]^{c_k},
$$
for $c_k=\sum_{i=1}^{j} 1\{n_i=k\}$, $k=1, \dots, n$. The {\it $[0:k]$-discovery} would therefore simply follow as the conditional probability:
$$
\hat{U}_{n+0}(k)=c_k \frac{q(c_1, \dots, c_k-1, c_{k+1}+1, \dots, c_{n+1})}{q(c_1, \dots, c_k, \dots, c_n)}= 
$$
\begin{equation}
\label{uncondo}=c_k \frac{V_{n+1,j}}{V_{n,j}}\frac{[(1-\alpha)_k]}{[(1-\alpha)_{k-1}]}=c_k\frac{V_{n+1,j}}{V_{n,j}} (k-\alpha).
\end{equation}

In what follows, as in Lijoi, Pr\"unster and Walker (2008) and Cerquetti (2013), given the allocation of the first $n$ observations in $j$ species with multiplicities $(n_1, \dots, n_j)$, let $K_m^{(n)}$ be the random number of {\it new} species arising in the additional $m$ sample, $(S_1, \dots, S_{K_m^{(n)}})$ the random vector of the multiplicities of the {\it new} species in {\it exchangeable random order} (Pitman, 2006, Eq. 2.7), $S_{m}=\sum_{i=1}^{K_m^{(n)}} S_{i}$ the total number of {\it new} observations belonging to {\it new} species and $(M_{1,m}, \dots, M_{j,m})$ the random vector of the allocation of the additional observations in the $j$ {\it old} species for $\sum_{i=1}^{j} M_{i,m}=m -S_m$.

\section{Correct Estimation of the $m$-Step $k$-Discovery}
\label{s:model}

\subsection{Under General Gibbs priors}

From (\ref{uncondo}) notice that, given a basic sample $(n_1, \dots, n_j)$, but assuming an intermediate $m$-sample still to be observed, the probability to observe at step $n+m+1$ a species represented $k$ times among the {\it new} species  generated by the additional sample, can be expressed as the random quantity
\begin{equation}
\label{newelle}
P^{n+m+1}_{new,k}(\alpha, V)=W_{k, m}^{(n)}\frac{V_{n+m+1, j+K_m^{(n)}}}{V_{n+m, j +K_m^{(n)}}} (k- \alpha),
\end{equation}
for $K_m^{(n)}$ the random number of {\it new} species induced by the additional sample and $W_{k,m}^{(n)}= \sum_{i=1}^{K_m^{(n)}}1 \{S_i=k|K_n=j\}$ the random number of {\it new} species represented $k$ times.

The next result has been established in Cerquetti (2013), cf. Propositions 8, by means of a far more easy proof of the one proposed in Favaro et al. (2012). For the sake of clarity we recall that {\it non central generalized Stirling numbers} and {\it non central generalized factorial coefficients} used in Favaro et al. (2012) are tied by the relation
$S_{n,\xi}^{-1, -\alpha, -\beta}= \mathcal{C}(n, \xi;\alpha,  -\beta)/\alpha^{\xi}$. We even recall the compact notation for generalized rising factorials $(x)_{n \uparrow y}=(x)(x+y)(x+2y)\cdots(x+(n-1)y)$. Notice also that in formula (\ref{newdisco}) may be useful to write the summation limits in the unconventional way $\sum_{l-1=0}^{m-k}$, as it is ubiquitous in Cerquetti (2013), to highlight that the correspondence with the parameters of the generalized Stirling numbers must be satisfied for the summation to make sense, and to exploit equation (\ref{stircon}) when needed. 

\bigskip

\noindent{\bf Theorem 1.} (Correction to Theorem 2, FLP 2012){\it  Under a general $(\alpha, V)$ Gibbs prior, for $W_{k,m}^{(n)}= \sum_{i=1}^{K_m^{(n)}}1 \{S_i=k|K_n=j\}$, the Bayesian nonparametric estimator of $P^{m+n+1}_{new,k}(\alpha, V)$ is given by 
\begin{equation}
\label{newdisco}
\hat{U}^{new}_{n+m}(k)=
(1 -\alpha)_k {m \choose k} \sum_{l=1}^{m-k+1} \frac{V_{n+m+1, j+l}}{V_{n,j}}   S_{m-k, l-1}^{-1, -\alpha, -(n-j\alpha)},
\end{equation}
where $S_{n,\xi}^{-1, -\alpha, -\gamma}$ are non central generalized Stirling numbers defined as the connection coefficients (cfr. e.g. Hsu and Shiue, 1998)
\begin{equation}
\label{stircon}
(x-\gamma)_n= \sum_{\xi=0}^n S_{n,\xi}^{-1, -\alpha, \gamma} (x)_{\xi \uparrow \alpha}
\end{equation}
for $(x)_{\xi \uparrow \alpha}= x(x+\alpha)\cdots(x +(\xi-1)\alpha)$ generalized rising factorials.}

\bigskip
\noindent{\it Proof}. See Cerquetti (2013), proof of Proposition 8.

\bigskip

\noindent{\bf Remark 2.} The problem in Theorem 2 in Favaro et al. (2012) is in the summation limits. The marginalization over the possible number $l$ of {\it new} species arising in the additional $m-$sample, giving rise to at least one species represented $k$ times,  must range, as in formula (\ref{newdisco}), between $1$ (for $m=k)$ and $m-k+1$ (for $m > k)$, since we can actually observe one species of size $k$ and $m-k$ species of size $1$. In Favaro et al. (2012) the second term of formula (12), written explicitly,  corresponds instead to 
$$
(1-\alpha)_k {m \choose k} \sum_{l=1}^{m-k} \frac{V_{n+m+1, j+l}}{V_{n,j}} \frac{1}{\alpha^{l-1}} \mathcal{C}(m-k, l-1; \alpha, -(n-j\alpha))
$$
where $\mathcal{C}(n,\xi; \alpha, -\beta)$ stands for {\it generalized factorial coefficients} (Chalarambides, 2005),
showing the summation ranges until $m-k$, thus producing the wrong result. 

\bigskip

\noindent{\bf Remark 3.} For $c_{i}$ the number of species observed $i$ times in the basic $n-$sample, $O_{k,m}^{(n)}= \sum_{i=1}^{j} 1\{n_i+M_{i,m}=k|n_1, \dots, n_j\}$ and $\mathbb{P}(M_{1,m}=m_{1}, \dots, M_{j,m}=m_j, S_m=s|n_1, \dots, n_j)$ {\it Multivariate Polya-Gibbs distributions} (cfr. Cerquetti, 2013, Sect. 4) arising by the conditional allocation of {\it new} observations among  {\it old} species, details of a far more easy derivation of the one proposed in Favaro et al. (2012), for the Bayesian nonparametric estimator, under general Gibbs prior, of the probability to observe at step $n+m+1$ a species represented $k$ times among  the {\it old} species, can be found in Cerquetti (2013, cf. Proposition 9). Here we report the result.  
\begin{equation}
\label{discoold}
\hat{U}_{n+m}^{old}(k)
= \sum_{i=1}^{k} c_i {m \choose {k-i}} (i-\alpha)_{k-i+1}
\times \sum_{l=0}^{m-k+i} \frac{V_{n+m+1, j+l}}{V_{n,j}} S_{m-k+i, l}^{-1, -\alpha, -(n-j\alpha - i +\alpha)},
\end{equation}
thus confirming the first term in Equation (12) in Favaro et al. (2012) is correct.

\subsection{Under Two-Parameter Poisson-Dirichlet Priors}

The problem in the general result under Gibbs priors affects the correctness of the result in Proposition 3, where the estimator under two-parameter Poisson Dirichlet priors is derived. The following correction is needed. Even in this case we are able to provide a considerably more easy proof with respect to the one in Favaro et al. (2012, cfr. Web Appendix pag. 8).

\bigskip

\noindent{\bf Proposition 4.} (Correction to Proposition 3. FLP 2012) {\it Specializing equation (\ref{newdisco}) and (\ref{discoold}) for the two-parameter Poisson-Dirichlet $(\alpha, \theta)$ model yields 
$$
\hat{U}_{n+m}^{\alpha, \theta}(k)= 
\sum_{i=1}^{k} c_i {m \choose k-i} \frac{(\theta +n -i +\alpha)_{m-k+i} (i-\alpha)_{k+1-i}}{(\theta +n)_{m+1}}+
$$
\begin{equation}
\label{poidirest}+\frac{(\theta +j\alpha)}{(\theta +n)} {m \choose k}  \frac{(1-\alpha)_{k} (\theta +\alpha +n)_{m-k}}{(\theta +n +1)_m}.
\end{equation}
}

\bigskip

\noindent{\it Proof}. By the specific form of the $(V_{n,j})$ coefficients for the two-parameter $(\alpha, \theta)$ model, for $\alpha \in (0,1)$ and $\theta > -\alpha$,
$$V_{n,j}^{(\alpha, \theta)}=\frac{(\theta +\alpha)_{j-1 \uparrow \alpha}}{(\theta +1)_{n-1}},$$
and the multiplicative property of generalized rising factorials, $(x)_{a+b \uparrow c}=(x)_{a \uparrow c} (x+ac)_{b \uparrow c}$,  equation (\ref{newdisco}) yields
$$
\hat{U}_{n+m}^{(\alpha, \theta), new}=\frac{(1-\alpha)_{k}}{(\theta +n)_{m+1}} {m \choose k} \sum_{l-1=0}^{m-k} (\theta +j\alpha)_{l \uparrow \alpha} S_{m-k, l-1}^{-1, -\alpha, -(n-j\alpha)}=
$$
and by the same property, and the definition of non-central generalized Stirling numbers as connection coefficients (\ref{stircon}),
$$
=\frac{(\theta +j\alpha)}{(\theta +n)} {m \choose k}  \frac{(1-\alpha)_{k} (\theta +\alpha +n)_{m-k}}{(\theta +n +1)_m}.
$$
Again, by the specific form of the $(\alpha, \theta)$ Gibbs coefficients, equation (\ref{discoold}) yields
$$ 
\hat{U}_{n+m}^{(\alpha, \theta), old}=\sum_{i=1}^k c_i \frac{(i-\alpha)_{k+1-i}}{(\theta +n)_{m+1}} {m \choose k-i}  \sum_{l=0}^{m-k+i} (\theta +j\alpha)_{l \uparrow \alpha} S_{m-k+i, l}^{-1, -\alpha, -(n -j\alpha -i +\alpha)}=
$$
and by the definition of non-central generalized Stirling numbers as connection coefficients (\ref{stircon})
$$
=\sum_{i=1}^{k} c_i {m \choose k-i} \frac{(\theta +n -i +\alpha)_{m-k+i} (i-\alpha)_{k+1-i}}{(\theta +n)_{m+1}}.
$$

\bigskip

\noindent{\bf Remark 5.} Relying on the proposed correction it turns out that the result in Favaro et al. (2012) underestimates the probability of interest under general Gibbs priors of an amount equal to 
$$
(1 -\alpha)_{k} {m \choose k} \frac{V_{n+m+1, j+m-k+1}}{V_{n,j}} S_{m-k, m-k}^{-1, -\alpha, -(n-j\alpha)},
$$
that reduces to
$$
 (1 -\alpha)_{k} {m \choose k} \frac{V_{n+m+1, j+m-k+1}}{V_{n,j}},
$$
since $S_{n,n}^{-1, -\alpha, -\gamma}=S_{n,n}^{-1, -\alpha}=1$. Under $(\alpha, \theta)$ Poisson-Dirichlet priors this quantity  corresponds to 
$$
{m \choose k}  \frac{(1-\alpha)_{k}}{(\theta +n)_{m+1}} (\theta +j\alpha)_{m-k+1 \uparrow \alpha},
$$
which is easily seen to be the term subtracted in Equation (16) in Proposition 3. in Favaro et al. (2012).

\subsection{A Counterproof}

To check the correctness of our correction, we verify that 
$$
\hat{U}_{n+m}^{\alpha, \theta}(0)+\sum_{k=1}^{m+n} \hat{U}_{n+m}^{\alpha, \theta}(k)= 1
$$
where $\hat{U}_{n+m}^{\alpha, \theta}(0)$ is the Bayesian nonparametric estimator of the  "discovery probability", the probability to discover a {\it new} species at observation $X_{n+m+1}$ under two-parameter Poisson-Dirichlet prior, as already obtained in Favaro, Lijoi, and Pr\"unster (2009, cfr. Eq. (7))
$$\hat{U}_{n+m}^{\alpha, \theta}(0)=\frac{(\theta +j\alpha)(\theta +\alpha +n)_m}{(\theta +n)(\theta+n+1)_m}.$$

\bigskip

\noindent{\bf Proposition 6.}
{\it Under two-parameter Poisson-Dirichlet $(\alpha, \theta)$ priors}
$$
\sum_{k=1}^{m+n} \hat{U}_{n+m}^{\alpha, \theta}(k)= 1 - \frac{(\theta +j\alpha)}{(\theta +n)} \frac{(\theta +\alpha +n)_m}{(\theta +n +1)_m}.
$$

\bigskip

\noindent{\it Proof}. With respect to the components given by the {\it new} species, possible sizes range between $1$ and $m$, hence, by the second term in (\ref{poidirest}), we calculate
$$
\frac{(\theta +j \alpha)}{(\theta +n)} \sum_{k=1}^{m}  {m \choose k} \frac{( 1- \alpha)_k (\theta + \alpha +n)_{m -k}}{(\theta +n +1)_m}=
$$
and by the definition  of Beta-Binomial distribution  of parameters $(m, (1-\alpha), (\theta + \alpha +n))$ (see e.g. Johnson and Kotz, 2005)
$$
=\frac{(\theta +j\alpha)}{(\theta +n)} \left( 1 - \frac{(\theta + \alpha +n)_m}{(\theta +n +1)_m}\right) =
$$
$$= \frac{(\theta +j\alpha)}{(\theta +n)} - \frac{(\theta +j\alpha)}{(\theta +n)} \frac{(\theta +\alpha +n)_m}{(\theta +n +1)_m}.
$$
It remains to check that summing the first term in (\ref{poidirest}) over possible sizes of the {\it old} species yields ${(n-j\alpha)}/{(\theta +n)}$, and in fact
$$
\sum_{k=1}^{n+m} \sum_{i=1}^{k} c_i{m \choose k-i} \frac{(\theta +n -i +\alpha)_{m-k+i}  (i-\alpha)_{k+1-i}}{(\theta +n)_{m+1}}=
$$
$$
=\sum_{i=1}^n  c_i  \sum_{k-i=0}^{m} {m \choose k-i} \frac{(\theta +n -i +\alpha)_{m-k+i}  (i-\alpha)_{k+1-i}}{(\theta +n)_{m+1}}=
$$
and by the multiplicative property of rising factorials and  the definition of Beta-Binomial distribution of parameters $(m, (i-\alpha+1), (\theta +n-i +\alpha))$ 
$$
=\sum_{i=1}^n  c_i \frac{(i-\alpha)}{(\theta +n)}
\times$$
$$\times  \sum_{k-i=0}^{m} {m \choose k-i} \frac{(\theta +n -i +\alpha)_{m-k+i}  (i-\alpha+1)_{k-i}}{(\theta +n+1)_{m}}=
$$
$$= \sum_{i=1}^n  c_i \frac{(i-\alpha)}{(\theta +n)}= \frac{(n-j\alpha)}{(\theta +n)}.
$$

\subsection {Correction to the Matlab Code}

From our correction the Matlab code provided by Favaro et al. (2012) in the Web Appendix (cfr. folder \verb+BIOM_1793_sm_suppmat.zip+, file {\sf discovery.m}) to estimate the {\it $[m:k]$ discovery} under $(\alpha, \theta)$ Poisson-Dirichlet priors should read

\begin{verbatim}BNP_mk_discovery_estimate(k,i)=
partial_sum+(exp(rising_factorial(1-sigma,k))
*binomial_coefficient(i,k)*(((theta+(j*sigma))
*exp(rising_factorial(theta+n+sigma,i-k)-
rising_factorial(theta+n,i+1)))));
\end{verbatim}
Estimates obtained in Table 3. in Section 3. ''An Application to Genomic Data'' should therefore be recalculated accordingly.

\section*{References}
\newcommand{\bibu}{\item \hskip-1.0cm}
\begin{list}{\ }{\setlength\leftmargin{1.0cm}}

\bibu \textsc{Cerquetti, A.} (2013). Marginals of multivariate Gibbs distributions with applications in Bayesian species sampling. {\it Electron. J. Statist.}, 7, 697-716. 

\bibu\textsc{Charalambides, C. A.} (2005). {\it Combinatorial Methods in Discrete Distributions}. Wiley, Hoboken NJ.

\bibu \textsc{Favaro, S., Lijoi, A., Mena, R. H., and Pr\"unster, I.} (2009). Bayesian non-parametric inference for species variety with a two-parameter Poisson-Dirichlet process prior. {\it Journal of the Royal Statistical Society B}, \textbf{71}, 993--1008.

\bibu\textsc{Favaro, S., Lijoi, A., and Pr\"unster, I.} (2012). A new estimator of the discovery probability. {\it Biometrics}, {\bf 68}, 4, 1188--1196. 

\bibu\textsc{Gnedin, A. and Pitman, J. } (2006). {Exchangeable Gibbs partitions  and Stirling triangles.} {\it Journal of Mathematical Sciences}, {\bf 138}, 3, 5674--5685. 

\bibu \textsc{Hsu, L. C, and Shiue, P. J.} (1998). A unified approach to generalized Stirling numbers. {\it Advances in Applied Mathematics}, {\bf 20}, 366-384.

\bibu \textsc{Johnson, N. S. and Kotz, S.} (2005). {\it Univariate discrete distributions} 3rd Ed. Wiley, NY. 

\bibu \textsc{Kingman, J.F.C} (1978). The representation of partition structure.  {\it Journal of London Mathematical Society} {\bf 2}, 374--380.

\bibu \textsc{Lijoi, A., Mena, R. H., and Pr\"unster, I.} (2007). Bayesian nonparametric estimation of the probability of discovering new species.  {\it Biometrika}, {\bf 94}, 769--786.

\bibu\textsc{Lijoi, A., Pr\"unster, I., and Walker, S.G.} (2008). Bayesian nonparametric estimators derived from conditional Gibbs structures. {\it Annals of Applied Probability}, \textbf{18}, 1519--154

\bibu\textsc{Pitman, J.} (2006). {\it Combinatorial Stochastic Processes}. Ecole d'Et\'e de Probabilit\'e de Saint-Flour XXXII - 2002. Lecture Notes in Mathematics N. 1875, Springer.

\bibu\textsc{Pitman, J. and Yor, M.} (1997). The two-parameter Poisson-Dirichlet distribution derived from a stable subordinator. {\it Annals of Probability}, {\bf25}, 855--900.


\end{list}

\label{lastpage}

\end{document}